\documentclass[leqno,a4paper]{article}
\usepackage{amssymb,amsmath,amsthm,verbatim}

\title{\textbf{A new look at the classification }\\
\textbf{of the tri-covectors of a} $6$\textbf{-dimensional
symplectic space }}
\author{\textsc{J. Mu\~{n}oz Masqu\'e, L. M. Pozo Coronado}}
\date{}

\newtheorem{theorem}{Theorem}[section]

\newtheorem{lemma}[theorem]{Lemma}
\newtheorem{corollary}[theorem]{Corollary}

\newtheorem{definition}[theorem]{Definition}

\begin{document}

\maketitle

\begin{abstract}
Let $\mathbb{F}$ be a field of characteristic $\neq 2$ and $3$,
let $V$ be a $\mathbb{F}$-vector space of dimension $6$,
and let $\Omega \in \wedge ^2V^\ast $ be a non-degenerate
form. A system of generators for polynomial invariant functions
under the tensorial action of the group $Sp(\Omega )$ on
$\wedge ^3 V^\ast $, is given explicitly. Applications
of these results to the normal forms of De Bruyn-Kwiatkowski
and Popov are given.
\end{abstract}

\bigskip


\medskip


\section{Symplectic invariants\label{symplectic}}
Below, $\mathbb{F}$ denotes a field of characteristic $\neq 2$ and $3$,
and $V$ is a $\mathbb{F}$-vector space of dimension $6$. Notations
and elementary properties of algebraic sets and groups have been
taken from Fogarty's book \cite{Fogarty}.

The group $GL(V)$ acts on $\otimes ^rV^\ast $ by
\[
(A\cdot \xi )(x_1,\dotsc,x_r)=\xi (A^{-1}x_1,\dotsc,A^{-1}x_r),
\]
$\forall\xi  \in \otimes ^rV^\ast ,\;\forall x_1,\dotsc,x_r\in V$,
and $GL(V)$ acts on $\wedge ^rV^\ast $ by the same formula.

In particular, the $\mathbb{F}$-homomorphism induced by the action
of $GL(V)$ on $\wedge ^3 V^\ast $ is denoted by
$\rho\colon GL(V)\to GL(\wedge ^3V^\ast )$.

We denote by $\rho^{\prime}\colon Sp(\Omega )
\to GL(\wedge ^3 V^\ast )$ the restriction of $\rho $
to the symplectic group of a non-degenerate
$2$-covector $\Omega \in \wedge ^2V^\ast $. Furthermore, $GL(V)$
acts on $\otimes ^rV^\ast \otimes V$ by
$(A\cdot \eta )(x_1,\dotsc,x_r)
=A[\eta (A^{-1}x_1,\dotsc,A^{-1}x_r)]$, for all $\eta $
in $\otimes ^rV^\ast \otimes V$, and all $x_1,\dotsc,x_r\in V $.
If $(v_i)_{i=1}^6$ is a basis of $V$ such that
$\Omega =\sum _{i=1}^6v^i\wedge v^{i+3}$ and $(v^i)_{i=1}^6$
is the dual basis, then we define a system
of coordinate functions $y_{abc}$,
$1\leq a<b<c\leq6$, on $\wedge ^3 V^\ast $ by setting
\begin{equation*}
\theta =\sum\nolimits_{1\leq a<b<c\leq 6}y_{abc}(\theta )
\left( v^a\wedge v^b\wedge v^c\right)
\in \wedge ^3 V^\ast .
\end{equation*}
For every $A\in GL(V)$ and every $1\leq a<b<c\leq 6$ we have

\begin{align}
A\cdot \left( v^a\wedge v^b\wedge v^c\right)
& =(A^{-1})^\ast v^a\wedge (A^{-1})^\ast v^b
\wedge (A^{-1})^\ast v^c \nonumber \\
& =\left( v^a\circ A^{-1}\right)
\wedge \left( v^b\circ A^{-1}\right)
\wedge \left( v^{c}\circ A^{-1}\right) \nonumber \\
& =\left( \lambda _{ah}v^h\right)
\wedge \left( \lambda _{bi}v^i\right)
\wedge \left( \lambda _{cj}v^{j}\right) \nonumber \\
& =\sum\nolimits_{1\leq h<i<j\leq 6}\left\vert
\begin{array}
[c]{ccc}
\lambda _{ah} & \lambda _{bh} & \lambda _{ch}\\
\lambda _{ai} & \lambda _{bi} & \lambda _{bi}\\
\lambda _{aj} & \lambda _{bj} & \lambda _{cj}
\end{array}
\right\vert v^h\wedge v^i\wedge v^j, \nonumber
\end{align}
where $(\lambda _{ij})_{i,j=1}^6$ is the matrix
of $(A^{-1})^T$ in the basis $(v_i)_{i=1}^6$. Therefore,
we have $\mathbb{F}[\wedge ^3 V^\ast ]
=\mathbb{F}[y_{abc}]_{1\leq a<b<c\leq 6}$, and hence,
$\mathbb{F}(\wedge ^3V^\ast )
=\mathbb{F}(y_{abc})_{1\leq a<b<c\leq 6}$.

A function $I\in\mathbb{F}[\wedge ^3 V^\ast ]$
(resp.\ $I\in \mathbb{F}(\wedge ^3 V^\ast )$)
is $Sp(\Omega )$-invariant if
\[
I\left( A\cdot \theta \right)
=I(\theta ),\quad \forall \theta \in \wedge ^3 V^\ast,
\; \forall A\in Sp(\Omega ).
\]
\section{The basic invariants defined\label{invariants}}
\subsection{$I_1$ defined\label{I_1}}
For every $\theta \in \wedge ^3 V^\ast $ there exists
a unique $J^\theta \in\wedge ^2V^\ast \otimes V$ such that
\begin{equation}
\theta (x,y,z)=\Omega (J^\theta (x,y),z),
\quad \forall x,y,z\in V.\label{J}
\end{equation}
Given $A\in GL(V)$ and replacing $\theta $ by $A\cdot \theta $
in the formula \eqref{J}, we have
\[
(A\cdot \theta )(x,y,z)=\Omega (J^{A\cdot \theta }(x,y),z).
\]
Expanding on the right-hand side, we deduce
\begin{align*}
(A\cdot \theta )(x,y,z) & =\theta (A^{-1}x,A^{-1}y,A^{-1}z)\\
& =\Omega (J^\theta (A^{-1}x,A^{-1}y),A^{-1}z)\\
& =(A\cdot \Omega)(A[J^\theta (A^{-1}x,A^{-1}y)],z))\\
& =(A\cdot \Omega)((A\cdot J^\theta )(x,y),z).
\end{align*}
Furthermore, if $A\in Sp(\Omega )$, then $A\cdot \Omega =\Omega $
and consequently
\[
(A\cdot \theta )(x,y,z)=\Omega ((A\cdot J^\theta )(x,y),z)
=\Omega (J^{A\cdot \theta }(x,y),z),
\quad \forall x,y,z\in V.
\]
Hence $J^{A\cdot \theta }=A\cdot J^\theta ,
\quad \forall A\in Sp(\Omega )$. If
\begin{equation}
\theta =\sum _{1\leq i<j<k\leq 6}
\lambda _{ijk}v^i\wedge v^{j}\wedge v^k
\label{theta}
\end{equation}
and $J^\theta =\sum_{b<c}\mu _{bc}^av^b\wedge v^c\otimes v_a$,
then by letting $x=v_i$, $y=v_j$, $z=z^kv_k$ in the formula
\eqref{J} for every pair $1\leq i<j\leq 5$, it follows:
\[
\sum _{k\neq i,j}\theta \left( v_i,v_j,v_k\right) z^k
=-\mu_{ij}^4z^1-\mu_{ij}^5z^2-\mu_{ij}^6z^3+\mu_{ij}^1z^4
+\mu_{ij}^2z^5+\mu_{ij}^3z^6,
\]
and comparing the coefficients of $z^1,\dotsc,z^6$ in both sides,
we deduce

\[
\begin{array}
[c]{lll}
\mu_{ij}^1=\lambda _{ij4},
& \mu_{ij}^2=\lambda _{ij5},
& \mu_{ij}^3=\lambda _{ij6}, \medskip\\
\mu_{ij}^4=-\lambda _{ij1},
& \mu_{ij}^{5}=-\lambda _{ij2},
& \mu_{ij}^6=-\lambda _{ij3},
\end{array}
\]
with the usual agreement:
$\lambda _{ijk}=\varepsilon _\sigma \lambda _{abc}$,
where $a<b<c$, $\{ a,b,c\} =\{ i,j,k\} $, $\sigma $ being
the permutation $a\mapsto i$, $b\mapsto j$, $c\mapsto k$.

Let $\Omega =\sum _{i=1}^6v^i\wedge v^{i+3}$ be
as in section \ref{symplectic}. Each
$\theta \in \wedge ^3 V^\ast $ determines a vector
$v_\theta \in V$ defined by the following equation:
\begin{equation}
i_{v_\theta }(\Omega \wedge \Omega \wedge \Omega )
=\theta\wedge \Omega .\label{v_theta}
\end{equation}
Transforming the equation \eqref{v_theta} by
$A\in Sp(\Omega )$ and recalling that
$A\cdot \Omega =\Omega $, we obtain
$A\cdot \left[
i_{v_\theta }(\Omega \wedge \Omega \wedge \Omega )
\right]
=(A\cdot \theta )\wedge \Omega
=i_{v_{A\cdot \theta }}(\Omega \wedge \Omega \wedge \Omega )$.
Moreover, for every system $x_1,\dotsc,x_5\in V$ one has:
\[
\begin{array}
[c]{rl}
\left(
A\cdot
\left[ i_{v_\theta }(\Omega \wedge \Omega \wedge \Omega )
\right]
\right)
\left( x_1,\dotsc,x_5\right)
=\! & \!
\left[ i_{v_\theta }(\Omega \wedge \Omega \wedge \Omega)
\right] \left( A^{-1}x_1,\dotsc,A^{-1}x_5\right) \\
=\! & \!(\Omega \wedge \Omega \wedge \Omega )
\left( A^{-1}Av_\theta ,A^{-1}
x_1,\dotsc,A^{-1}x_5\right) \\
=\! & \!\!
\left( i_{Av_\theta }\left[
A\cdot (\Omega \wedge \Omega \wedge \Omega )\right]
\right) \left( x_1,\dotsc,x_5\right) .
\end{array}
\]
Hence $A\cdot
\left[
i_{v_\theta }(\Omega \wedge \Omega \wedge \Omega )
\right]
=i_{Av_\theta }
\left[
A\cdot (\Omega \wedge \Omega \wedge \Omega )
\right]
=i_{Av_\theta }(\Omega \wedge \Omega\wedge \Omega)$.
Accordingly: $v_{A\cdot \theta }=Av_\theta $.
Letting $x=v_\theta $ in \eqref{J}, it follows:
\[
\Omega \left[ \left( i_{v_\theta }J^\theta \right) (y),z\right]
=\Omega\left[ J^\theta (v_\theta ,y),z\right]
=\theta (v_\theta ,y,z),
\]
and replacing $\theta $ by $A\cdot \theta $, $A\in Sp(\Omega )$,
we have

\begin{align*}
\Omega
\left[
\left(
i_{v_{A\cdot \theta }}J^{A\cdot \theta }
\right) (y),z
\right]
& =\Omega
\left[
J^{A\cdot \theta }(v_{A\cdot \theta },y),z
\right]
=(A\cdot \theta )(v_{A\cdot \theta },y,z)\\
& =\theta
\left(
A^{-1}v_{A\cdot \theta },A^{-1}y,A^{-1}z
\right)
=\theta (v_\theta ,A^{-1}y,A^{-1}z)\\
& =\Omega
\left[
J^\theta (v_\theta ,A^{-1}y),A^{-1}z
\right]
=\Omega
\left[
(i_{v_\theta }J^\theta )(A^{-1}y),A^{-1}z
\right] \\
& =\left( A\cdot \Omega
\right)
\left[
(i_{v_\theta }J^\theta )(y),z
\right]
\\
& =\Omega
\left[
(i_{v_\theta }J^\theta )(y),z
\right] .
\end{align*}
Hence $i_{v_{A\cdot \theta }}J^{A\cdot \theta }
=i_{v_\theta }J^\theta $ for all $A\in Sp(\Omega )$.
Accordingly, the endomorphism
$L_\theta =i_{v_\theta }J^\theta \in V^\ast \otimes V$
is $Sp(\Omega )$-invariant.

If $\theta $ is as in \eqref{theta} and
$v_\theta =\sum _{h=1}^6x^hv_h$, then,
as a computation shows, we have
\begin{equation*}
\begin{array}
[c]{lll}
x^1=-\tfrac{1}{6}(\lambda _{245}+\lambda _{346}),
& x^2=\tfrac{1}{6}(\lambda _{145}-\lambda _{356}),
& x^3=\tfrac{1}{6}(\lambda _{146}+\lambda _{256}),
\medskip\\
x^4=-\tfrac{1}{6}(\lambda _{125}+\lambda _{136}),
& x^5=-\tfrac{1}{6}(\lambda _{236}-\lambda _{124}),
& x^6=\tfrac{1}{6}(\lambda _{134}+\lambda _{235}).
\end{array}
\end{equation*}

Again writing
$J^\theta =\sum _{b<c}\mu _{bc}^av^b\wedge v^c\otimes v_a$,
we obtain
\begin{equation*}
\begin{array}
[c]{lll}
\mu_{ij}^1=\lambda _{ij4},
& \mu _{ij}^2=\lambda _{ij5},
& \mu _{ij}^3=\lambda _{ij6},\medskip\\
\mu _{ij}^4=-\lambda _{ij1},
& \mu _{ij}^{5}=-\lambda _{ij2},
& \mu_{ij}^6=-\lambda _{ij3},
\end{array}
\end{equation*}
or equivalently $\mu_{ij}^{h}=\lambda _{ij,h+1}$,
$\mu_{ij}^{h+3}=-\lambda _{ij,h}$ for $1\leq h\leq3$.
Therefore
\[
L_\theta =i_{v_\theta }J^\theta
=\sum\nolimits_{b<c}\mu_{bc}^a\left(
x^bv^c-x^cv^b\right) \otimes v_a.
\]

Hence, the matrix $M=(M_{ia})_{i,a=1}^6$ of $L_\theta $
in the basis $v_1,\dotsc,v_6$ is given by
$M_{ia}=\sum\nolimits_{b=1}^6\mu_{bi}^ax^b$,
and as a computation shows, the characteristic polynomial
of $L_\theta $ is written as follows:
\[
\det \left( xI-L_\theta \right)
=x^6+c_4(\theta)x^4+c_2(\theta )x^2,
\]
where
\[
\begin{array}
[c]{ll}
c_2(\theta )=\frac{1}{2^4\cdot 3^4}(I_1)^2,
& c_4(\theta )=-\frac{1}{2\cdot3^2}I_1,
\end{array}
\]
and the explicit expression of $I_1$ is given in the Appendix.
\subsection{$I_2$ defined\label{I_2}}
For every $\theta \in \wedge ^3 V^\ast $ let
$J^\theta \barwedge J^\theta \in\wedge ^3 V^\ast \otimes V$
be the tensor defined as follows:
\[
\begin{array}
[c]{ll}
\left(
J^\theta \barwedge J^\theta
\right) \left( x_1,x_2,x_3\right)
= & J^\theta \left( J^\theta \left( x_1,x_2\right) ,x_3\right)
+J^\theta
\left( J^\theta \left( x_2,x_3\right) ,x_1\right) \\
\multicolumn{1}{r}{} & \multicolumn{1}{r}{+J^\theta
\left( J^\theta \left( x_3,x_1\right) ,x_2\right) ,}\\
\multicolumn{1}{r}{} & \multicolumn{1}{r}{\forall x_1,x_2,x_3\in V.}
\end{array}
\]
As $J^\theta \left( v_i,v_j\right) =\mu_{ij}^hv_h$, we have
$J^\theta \left( J^\theta \left( v_i,v_j\right) ,v_k\right)
=\mu_{ij}^hJ^\theta \left( v_h,v_k\right)
=\mu_{ij}^h\mu _{hk}^lv_l$. Hence
\[
\left( J^\theta \barwedge J^\theta \right)
\left( v_i,v_j,v_k\right)
=\left( \mu_{ij}^h\mu_{hk}^l+\mu_{jk}^h\mu_{hi}^l
+\mu_{ki}^h\mu_{hj}^l\right) v_l.
\]

Let $\Omega ^\theta \in\wedge ^6V^\ast $ be defined by
$\Omega ^\theta =\tfrac{1}{3!^2}\operatorname{alt}
\left[ \tilde{\Omega }\circ \left(
J^\theta \barwedge J^\theta \otimes J^\theta \barwedge J^\theta
\right) \right] $, where
$\tilde{\Omega }\colon V\otimes V\to \mathbb{F}$ is
the linear map attached to $\Omega $; i.e.,
$\tilde{\Omega }(x\otimes y)=\Omega (x,y)$, $\forall x,y\in V$,
and $(J^\theta \barwedge J^\theta )
\otimes (J^\theta \barwedge J^\theta )
\colon \wedge ^3 V\otimes\wedge ^3V\to V\otimes V$
is the tensor product of $J^\theta \barwedge J^\theta $
and itself, $J^\theta \barwedge J^\theta $ being understood
as a linear map $\wedge ^3 V\to V$, via the canonical isomorphism
$\wedge ^3 V^\ast \otimes V=\operatorname{Hom}
\left( \wedge ^3 V,V\right) $.

Letting $\xi _{ijk}^l=\mu_{ij}^h\mu_{hk}^l+\mu_{jk}^h\mu_{hi}^l
+\mu_{ki}^h\mu_{hj}^l$, we have
\[
J^\theta \barwedge J^\theta
=\sum_{i<j<k}\xi _{ijk}^lv^i\wedge v^j\wedge v^k\otimes v_l.
\]
Hence
$\Omega ^\theta =\sum_{i<j<k}\sum_{a<b<c}\xi _{ijk}^l\xi _{abc}^d
\left( v^i\wedge v^j\wedge v^k\right) \wedge
\left( v^a\wedge v^b\wedge v^c\right) \Omega (v_l,v_d)$.

The product
$\left( v^i\wedge v^j\wedge v^k\right)
\wedge \left( v^a\wedge v^b\wedge v^c\right) $
does not vanish if and only if the indices $i,j,k,a,b,c$
are pairwise distinct, i.e., $\{ i,j,k,a,b,c\} =\{ 1,\dotsc,6\} $,
and in that case,
$v^i\wedge v^j\wedge v^k\wedge v^a\wedge v^b\wedge v^c
=\varepsilon _\sigma
v^1\wedge v^2\wedge v^3\wedge v^4\wedge v^5\wedge v^6$,
$\sigma $ being the permutation $\sigma (1)=i$,
$\sigma (2)=j$, $\sigma (3)=k$, $\sigma (4)=a$,
$\sigma (5)=b$, $\sigma (6)=c$. Therefore
$\Omega ^\theta =I_2(\theta)\Omega \wedge \Omega \wedge \Omega $,
where
\begin{align*}
I_2(\theta )\!& =\! -\tfrac{1}{6}\sum _{\sigma \in S_6}
\varepsilon _\sigma \xi _{\sigma (1)\sigma (2)\sigma (3)}^l
\xi _{\sigma (4)\sigma (5)\sigma (6)}^{d}\Omega (v_{l},v_{d}) \\
\!& =\! \tfrac{1}{6}\sum _{\sigma \in S_6}
\varepsilon _\sigma \xi _{\sigma (1)\sigma (2)\sigma (3)}^4
\xi _{\sigma (4)\sigma (5)\sigma (6)}^{1}
-\tfrac{1}{6}\sum _{\sigma \in S_6}
\varepsilon _\sigma \xi _{\sigma (1)\sigma (2)\sigma (3)}^1
\xi _{\sigma (4)\sigma (5)\sigma (6)}^4 \\
& \!\! +\tfrac{1}{6}\sum_{\sigma \in S_{6}}
\varepsilon _{\sigma }\xi _{\sigma (1)\sigma (2)\sigma (3)}^5
\xi _{\sigma (4)\sigma (5)\sigma (6)}^2
-\tfrac{1}{6}\sum_{\sigma \in S_6}
\varepsilon _{\sigma }\xi _{\sigma (1)\sigma (2)\sigma (3)}^2
\xi _{\sigma (4)\sigma (5)\sigma (6)}^5 \\
& \!\! +\tfrac{1}{6}\sum _{\sigma \in S_6}
\varepsilon _\sigma \xi _{\sigma (1)\sigma (2)\sigma (3)}^6
\xi _{\sigma (4)\sigma (5)\sigma (6)}^3
-\tfrac{1}{6}\sum_{\sigma \in S_{6}}
\varepsilon _\sigma \xi _{\sigma (1)\sigma (2)\sigma (3)}^3
\xi _{\sigma (4)\sigma (5)\sigma (6)}^6.
\end{align*}
If $\sigma _0$ is the permutation $1\mapsto 4$,
$2\mapsto 5$, $3\mapsto 6$, $4\mapsto 1$, $5\mapsto 2$,
$6\mapsto3$, then letting $\sigma ^\prime
=\sigma\circ\sigma_{0}$, we have $\sigma _0\circ \sigma _0=1$,
$\varepsilon _{\sigma ^\prime}
=\varepsilon _\sigma \varepsilon _{\sigma _0}
=-\varepsilon _\sigma $,

\[
\begin{array}
[c]{lll}
\sigma ^\prime (1)=\sigma (4), & \sigma ^\prime (2)
=\sigma (5), & \sigma ^\prime (3)=\sigma (6),\\
\sigma ^\prime (4)=\sigma (1), & \sigma ^\prime (5)
=\sigma (2), & \sigma ^\prime (6)=\sigma (3).
\end{array}
\]
and from the previous formula for $I_2(\theta )$ we obtain
\begin{equation*}
\begin{array}
[c]{rl}
I_2(\theta)= & -\tfrac{1}{3}\sum _{\sigma \in S_6}
\varepsilon _\sigma
\xi _{\sigma (1)\sigma (2)\sigma (3)}^1
\xi _{\sigma (4)\sigma (5)\sigma (6)}
^4\smallskip\\
& -\tfrac{1}{3}\sum_{\sigma\in S_6}
\varepsilon _\sigma
\xi _{\sigma (1)\sigma (2)\sigma (3)}^2
\xi _{\sigma (4)\sigma (5)\sigma (6)}^5
\smallskip\\
& -\tfrac{1}{3}\sum_{\sigma \in S_6}
\varepsilon _\sigma
\xi _{\sigma (1)\sigma (2)\sigma (3)}^3
\xi _{\sigma (4)\sigma (5)\sigma (6)}^6.
\end{array}
    \end{equation*}
The explicit expression for $I_2$ is also given in the Appendix.
\section{Infinitesimal criterion of invariance}
The derivations
\[
\begin{array}
[c]{l}
\tfrac{\partial }{\partial y_{abc}}\colon \mathbb{F}[\wedge ^3 V^\ast]
\to \mathbb{F}[\wedge ^3 V^\ast ]\quad(\mathrm{resp.}\;\tfrac
{\partial}{\partial y_{abc}}\colon\mathbb{F}(\wedge ^3 V^\ast )
\to \mathbb{F}(\wedge ^3 V^\ast )),\\
1\leq a<b<c\leq 6,
\end{array}
\]
are a basis of the $\mathbb{F}[\wedge ^3 V^\ast ]$-module (resp.\
$\mathbb{F}(\wedge ^3 V^\ast )$-vector space)
$\operatorname{Der}\nolimits_{\mathbb{F}}\mathbb{F}[\wedge ^3 V^\ast ]$
(resp.\ $\operatorname{Der}
\nolimits_{\mathbb{F}}\mathbb{F}(\wedge ^3 V^\ast )$). The subrings
of $Sp(\Omega )$-invariant functions are denoted by
$\mathbb{F}[\wedge ^3 V^\ast ]^{Sp(\Omega )}$
and $\mathbb{F}(\wedge ^3 V^\ast )^{Sp(\Omega )}$, respectively.
\begin{lemma}
\label{lemma_bis}
Assume $V$ is a $\mathbb{F}$-vector space of dimension $6$
and let $\Omega\in\wedge ^2V^\ast $ be a non-degenerate $2$-covector.
If
$I\in\mathbb{F}[\wedge ^3 V^\ast ]$ (resp.\
$I\in\mathbb{F}(\wedge ^3V^\ast )$) is a $Sp(\Omega )$-invariant
function, then $I$ is a common first integral of the following
derivations:
\begin{equation}
\begin{array}
[c]{l}
U^\ast =\sum\nolimits_{1\leq h<i<j\leq 6}
\left(
\sum\nolimits_{1\leq a<b<c\leq 6}U_{hij}^{abc}y_{abc}
\right)
\tfrac{\partial}{\partial y_{hij}},\\
\forall U=(u_{ij})_{i,j=1}^6\in \mathfrak{sp}(6,\mathbb{F}),
\end{array}
\label{U^ast}
\end{equation}
where the functions $U_{hij}^{abc}$ are defined by
\begin{equation*}
U_{hij}^{abc}=-\left\vert
\begin{array}
[c]{ccc}
u_{ha} & \delta_{hb} & \delta_{hc}\\
u_{ia} & \delta_{ib} & \delta_{ic}\\
u_{ja} & \delta_{jb} & \delta_{jc}
\end{array}
\right\vert -\left\vert
\begin{array}
[c]{ccc}
\delta_{ha} & u_{hb} & \delta_{hc}\\
\delta_{ia} & u_{ib} & \delta_{ic}\\
\delta_{ja} & u_{jb} & \delta_{jc}
\end{array}
\right\vert -\left\vert
\begin{array}
[c]{ccc}
\delta_{ha} & \delta_{hb} & u_{hc}\\
\delta_{ia} & \delta_{ib} & u_{ic}\\
\delta_{ja}^{{}} & \delta_{jb} & u_{jc}
\end{array}
\right\vert ,
\end{equation*}
and $\delta$ denotes the Kronecker delta.
\end{lemma}

\begin{proof}
We first observe that the derivations in \eqref{U^ast}
are the image of the $\mathbb{F}$-homomorphism
of Lie algebras $\rho _\ast ^\prime
\colon \mathfrak{sp}(6,\mathbb{F})\to \operatorname{Der}
\nolimits_{\mathbb{F}}\mathbb{F}[\wedge ^3 V^\ast ]$ induced
by $\rho ^\prime $. Accordingly, we only need to show that
$U^\ast (I)=0$ for the elements $U$ in a basis of the algebra
$\mathfrak{sp}(6,\mathbb{F})$. If $(E_{ij})_{i,j=1}^6$ is the
standard basis of $\mathfrak{gl}(6,\mathbb{F})$, then the matrices
\begin{equation}
\begin{array}
[c]{llll}
E_{11}+E_{14}-E_{41}-E_{44}, & E_{41}, & E_{52}, & E_{14},\\
E_{22}-E_{25}+E_{52}-E_{55}, & E_{42}+E_{51}, & E_{53}+E_{62}, & E_{25},\\
E_{33}-E_{66}+E_{36}-E_{63}, & E_{43}+E_{61}, & E_{63}, & E_{36},\\
E_{12}-E_{54}, & E_{31}-E_{46}, & E_{15}+E_{24}, & E_{26}+E_{35},\\
E_{21}-E_{45}, & E_{23}-E_{65}, & E_{16}+E_{34}, & E_{13}-E_{64},\\
E_{32}-E_{56}, &  &  &
\end{array}
\label{basis}
\end{equation}
are a basis $\mathcal{B}$ of $\mathfrak{sp}(6,\mathbb{F})$
with the following property: For every $U\in\mathcal{B}$,
we have $U^2=0$, as a simple computation proves. Hence,
for every $U\in\mathcal{B}$ and $t\in\mathbb{F}$,
the endomorphism $I+tU$ ($I$ denoting the identity map of $V$)
is symplectic. In fact, if $M_\Omega $, $M_U$ are the matrices
of $\Omega $, $U$, respectively, then
\begin{align*}
\left( M_{I+tU}\right) ^{T}M_{\Omega}M_{I+tU}
& =\left( I+t\left( M_U\right) ^T\right)
M_\Omega \left( I+tM_U\right) \\
& =M_\Omega +t\left[ \left( M_U\right) ^TM_\Omega
+M_\Omega M_U\right] +t^2\left( M_U\right) ^TM_\Omega M_U,
\end{align*}
but on one hand, we have
$\left( M_U\right) ^TM_\Omega +M_\Omega M_{U}=0$,
as $U\in\mathfrak{sp}(6,\mathbb{F})$, and on the other:
$\left( M_U\right) ^TM_\Omega M_U=-M_\Omega M_{U^2}=0$.
Hence $\left( M_{I+tU}\right) ^TM_\Omega M_{I+tU}
=M_\Omega $, thus proving that $I+tU\in Sp(\Omega )$,
$\forall t\in \mathbb{F}$.

Therefore $I\left( (I+tU)\cdot \theta \right) =I(\theta )$,
 $\forall t\in\mathbb{F}$,
 $\forall U=(u_{ij})\in \mathcal{B}$ and
 $\forall t\in \mathbb{F} $. If
 $\Lambda(t)=(\lambda _{ij}(t))_{i,j=1}^6=I-tU^T$,
then
\[
I\Bigl(
\sum\nolimits_{\substack{1\leq a<b<c\leq 6 \\1\leq h<i<j\leq 6}}~y_{abc}
\left\vert
\begin{array}
[c]{ccc}
\lambda _{ah}(t) & \lambda _{bh}(t) & \lambda _{ch}(t)\\
\lambda _{ai}(t) & \lambda _{bi}(t) & \lambda _{ci}(t)\\
\lambda _{aj}(t) & \lambda _{bj}(t) & \lambda _{cj}(t)
\end{array}
\right\vert v^{h}\wedge v^i\wedge v^{j}
\Bigr) =I(\theta ),
\]
and taking derivatives at $t=0$, we obtain
\[
0=\sum\nolimits_{1\leq a<b<c\leq 6,1\leq h<i<j\leq6}
U_{hij}^{abc}y_{abc}
\tfrac{\partial I}{\partial y_{hij}}(\theta).
\]
\end{proof}
\begin{definition}
Let $F$ be the field of fractions of an entire ring $R$.
The \emph{generic rank} of a finitely-generated $R$-module
$\mathcal{M}$ is the dimension of the $F$-vector space
$F\otimes_r\mathcal{M}$.
\end{definition}
\begin{theorem}
\label{theorem}
The generic rank of the $\mathbb{F}[\wedge ^3 V^\ast ]$-module
$\mathcal{M}$ spanned by the derivations in the formula
\emph{\eqref{U^ast}} of \emph{Lemma \ref{lemma_bis}} is $18$.
\end{theorem}
\begin{proof}
If $U=(u_{ij})_{i,j=1}^6\in\mathfrak{sp}(\Omega)$ then
\[
\begin{array}
[c]{lllll}
u_{24}=u_{15}, & u_{34}=u_{16},
& u_{35}=u_{26}, & u_{51}=u_{42}, &
u_{61}=u_{43},\\
u_{62}=u_{53}, & u_{44}=-u_{11},
& u_{45}=-u_{21}, & u_{46}=-u_{31}, &
u_{54}=-u_{12},\\
\multicolumn{1}{r}{u_{55}=-u_{22},}
& \multicolumn{1}{r}{u_{56}=-u_{32},}
& \multicolumn{1}{r}{u_{64}=-u_{13},}
& \multicolumn{1}{r}{u_{65}=-u_{23},} &
\multicolumn{1}{r}{u_{66}=-u_{33},}
\end{array}
\]
and the following $21$ functions are coordinates
on the symplectic algebra $\mathfrak{sp}(\Omega )$:
$u_{11}$, $u_{12}$, $u_{13}$, $u_{14}$, $u_{15}$,
$u_{16}$, $u_{21}$, $u_{22}$, $u_{23}$, $u_{25}$,
$u_{26}$, $u_{31}$, $u_{32}$, $u_{33}$, $u_{36}$,
$u_{41}$, $u_{42}$, $u_{43}$, $u_{52}$, $u_{53}$,
$u_{63}$. Hence we can write:
\[
\begin{array}
[c]{rl}
U^\ast \! = & \!\!\!\!
u_{11}Z_{11}+u_{12}Z_{12}+u_{13}Z_{13}+u_{14}Z_{14}
+u_{15}Z_{15}+u_{16}Z_{16}+u_{21}Z_{21}\\
\! & \multicolumn{1}{r}{\!\!\!\!
+u_{22}Z_{22}+u_{23}Z_{23}+u_{25}Z_{25}+u_{26}Z_{26}
+u_{31}Z_{31}+u_{32}Z_{32}+u_{33}Z_{33}}\\
\! & \multicolumn{1}{r}{\!\!\!\!+u_{36}Z_{36}
+u_{41}Z_{41}+u_{42}Z_{42}+u_{43}Z_{43}
+u_{52}Z_{52}+u_{53}Z_{53}+u_{63}Z_{63},}
\end{array}
\]
where
\begin{equation*}
\begin{array}{rc}
Z_{11}= &
-y_{123}Y_{123}-y_{125}Y_{125}-y_{126}Y_{126}
-y_{135}Y_{135}-y_{136}Y_{136}-y_{156}Y_{156}
\\
&
+y_{234}Y_{234}+y_{245}Y_{245}+y_{246}Y_{246}
+y_{345}Y_{345}+y_{346}Y_{346}+y_{456}Y_{456},
\\
Z_{12}= &
y_{124}Y_{125}-y_{234}Y_{134}+y_{134}Y_{135}
-y_{235}Y_{135}-y_{236}Y_{136}-y_{245}Y_{145}
\\
&
-y_{246}Y_{146}+y_{146}Y_{156}-y_{256}Y_{156}
+y_{234}Y_{235}+y_{246}Y_{256}+y_{346}Y_{356},
\\
Z_{13}= &
y_{234}Y_{124}+y_{235}Y_{125}+y_{124}Y_{126}
+y_{236}Y_{126}+y_{134}Y_{136}-y_{345}Y_{145}
\\
&
-y_{346}Y_{146}-y_{145}Y_{156}-y_{356}Y_{156}
+y_{234}Y_{236}-y_{245}Y_{256}-y_{345}Y_{356},
\\
Z_{14}= &
-y_{234}Y_{123}+y_{245}Y_{125}+y_{246}Y_{126}
+y_{345}Y_{135}+y_{346}Y_{136}-y_{456}Y_{156},
\\
Z_{15}= &
y_{134}Y_{123}-y_{235}Y_{123}-y_{245}Y_{124}
-y_{145}Y_{125}-y_{146}Y_{126}+y_{256}Y_{126}
\\
&
-y_{345}Y_{134}+y_{356}Y_{136}+y_{456}Y_{146}
+y_{345}Y_{235}+y_{346}Y_{236}-y_{456}Y_{256},
\end{array}
\end{equation*}
\begin{equation*}
\begin{array}{rc}
Z_{16}= &
-y_{124}Y_{123}-y_{236}Y_{123}-y_{246}Y_{124}
-y_{256}Y_{125}-y_{346}Y_{134}-y_{145}Y_{135}
\\
&
-y_{356}Y_{135}-y_{146}Y_{136}-y_{456}Y_{145}
-y_{245}Y_{235}-y_{246}Y_{236}-y_{456}Y_{356},
\\
Z_{21}= &
y_{125}Y_{124}+y_{135}Y_{134}+y_{156}Y_{146}
-y_{134}Y_{234}+y_{235}Y_{234}-y_{135}Y_{235}
\\
&
-y_{136}Y_{236}-y_{145}Y_{245}-y_{146}Y_{246}
+y_{256}Y_{246}-y_{156}Y_{256}+y_{356}Y_{346},
\\
Z_{22}= &
-y_{123}Y_{123}-y_{124}Y_{124}-y_{126}Y_{126}
+y_{135}Y_{135}+y_{145}Y_{145}+y_{156}Y_{156}
\\
&
-y_{234}Y_{234}-y_{236}Y_{236}-y_{246}Y_{246}
+y_{345}Y_{345}+y_{356}Y_{356}+y_{456}Y_{456},
\\
Z_{23}= &
-y_{134}Y_{124}-y_{135}Y_{125}+y_{125}Y_{126}
-y_{136}Y_{126}+y_{135}Y_{136}+y_{145}Y_{146}
\\
&
+y_{235}Y_{236}-y_{345}Y_{245}+y_{245}Y_{246}
-y_{346}Y_{246}-y_{356}Y_{256}+y_{345}Y_{346},
\\
Z_{25}= &
y_{135}Y_{123}+y_{145}Y_{124}-y_{156}Y_{126}
-y_{345}Y_{234}+y_{356}Y_{236}+y_{456}Y_{246},
\\
Z_{26}= &
-y_{125}Y_{123}+y_{136}Y_{123}+y_{146}Y_{124}
+y_{156}Y_{125}+y_{145}Y_{134}-y_{156}Y_{136}
\\
&
+y_{245}Y_{234}-y_{346}Y_{234}-y_{356}Y_{235}
-y_{256}Y_{236}-y_{456}Y_{245}+y_{456}Y_{346},
\\
Z_{31}= &
y_{126}Y_{124}+y_{136}Y_{134}-y_{156}Y_{145}
+y_{124}Y_{234}+y_{236}Y_{234}+y_{125}Y_{235}
\\
&
+y_{126}Y_{236}-y_{256}Y_{245}-y_{145}Y_{345}
-y_{356}Y_{345}-y_{146}Y_{346}-y_{156}Y_{356},
\end{array}
\end{equation*}
\[
\begin{array}
[c]{rl}
Z_{32}= & y_{126}Y_{125}-y_{124}Y_{134}-y_{125}Y_{135}
+y_{136}Y_{135}-y_{126}Y_{136}+y_{146}Y_{145}\\
& \multicolumn{1}{r}{+y_{236}Y_{235}+y_{246}Y_{245}-y_{245}Y_{345}
+y_{346}Y_{345}-y_{246}Y_{346}-y_{256}Y_{356},}\\
Z_{33}= & -y_{123}Y_{123}+y_{126}Y_{126}-y_{134}Y_{134}
-y_{135}Y_{135}+y_{146}Y_{146}+y_{156}Y_{156}\\
& \multicolumn{1}{r}{-y_{234}Y_{234}-y_{235}Y_{235}+y_{246}Y_{246}
+y_{256}Y_{256}-y_{345}Y_{345}+y_{456}Y_{456},}\\
Z_{36}= & -y_{126}Y_{123}+y_{146}Y_{134}+y_{156}Y_{135}+y_{246}Y_{234}
+y_{256}Y_{235}-y_{456}Y_{345},\\
Z_{41}= & -y_{123}Y_{234}+y_{125}Y_{245}+y_{126}Y_{246}+y_{135}Y_{345}
+y_{136}Y_{346}-y_{156}Y_{456},\\
Z_{42}= & y_{123}Y_{134}-y_{125}Y_{145}-y_{126}Y_{146}-y_{123}Y_{235}
-y_{124}Y_{245}+y_{126}Y_{256}\\
& \multicolumn{1}{r}{-y_{134}Y_{345}+y_{235}Y_{345}+y_{236}Y_{346}
+y_{136}Y_{356}+y_{146}Y_{456}-y_{256}Y_{456},}\\
Z_{43}= & -y_{123}Y_{124}-y_{135}Y_{145}-y_{136}Y_{146}-y_{123}Y_{236}
-y_{235}Y_{245}-y_{124}Y_{246}\\
& \multicolumn{1}{r}{-y_{236}Y_{246}-y_{125}Y_{256}-y_{134}Y_{346}
-y_{135}Y_{356}-y_{145}Y_{456}-y_{356}Y_{456},}\\
Z_{52}= & y_{123}Y_{135}+y_{124}Y_{145}-y_{126}Y_{156}-y_{234}Y_{345}
+y_{236}Y_{356}+y_{246}Y_{456},\\
Z_{53}= & -y_{123}Y_{125}+y_{123}Y_{136}+y_{134}Y_{145}+y_{124}Y_{146}
+y_{125}Y_{156}-y_{136}Y_{156}\\
& +y_{234}Y_{245}-y_{236}Y_{256}-y_{234}Y_{346}-y_{235}Y_{356}-y_{245}
Y_{456}+y_{346}Y_{456},\\
Z_{63}= & -y_{123}Y_{126}+y_{134}Y_{146}+y_{135}Y_{156}+y_{234}Y_{246}
+y_{235}Y_{256}-y_{345}Y_{456},
\end{array}
\]
where $Y_{abc}=\frac{\partial }{\partial y_{abc}}$,
$1\leq a<b<c\leq6$, is the standard basis of derivations.

Moreover, the invariant functions $I_1$ and $I_2$ are algebraically
independent. In fact, by using the formulas for $I_1$ and $I_2$ in the
Appendix, after a computation, it follows that the determinant
\[
\left\vert
\begin{array}
[c]{cc}
dI_1(Y_{123}) & dI_1(Y_{126})\\
dI_2(Y_{123}) & dI_2(Y_{126})
\end{array}
\right\vert
\]
at the point
\[
\begin{array}
[c]{ccccccc}
\lambda _{123}=1,
& \lambda _{124}=0,
& \lambda _{125}=1,
& \lambda _{126}=0,
& \lambda _{134}=1,
& \lambda _{135}=0,
& \lambda _{136}=0,\\
\lambda _{145}=0,
& \lambda _{146}=0,
& \lambda _{156}=0,
& \lambda _{234}=0,
& \lambda _{235}=0,
& \lambda _{236}=0,
& \lambda _{245}=0,\\
\lambda _{246}=0,
& \lambda _{256}=0,
& \lambda _{345}=0,
& \lambda _{346}=0, &
\lambda _{356}=0,
& \lambda _{456}=1, &
\end{array}
\]
takes the value $2^4\cdot 3^2$. As the differentials $dI_1$ and $dI_2$
vanish over all the vector fields $Z_{11},\dotsc,Z_{63}$, we deduce
that the generic rank of $\mathcal{M}$ is $\leq 18$.

Moreover, it is easy to obtain values of the variables $y_{abc}$,
$1\leq a<b<c\leq 6$, for which the matrix of $Z_{11},\dotsc,Z_{63}$
in the basis $Y_{abc}$, $1\leq a<b<c\leq6$, is $18$, and thus
we can conclude the proof.
\end{proof}
\begin{corollary}
If $\mathbb{F}$ is an algebraically closed field
of characteristic zero, then
$\mathbb{F}[\wedge ^3 V^\ast ]^{Sp(\Omega )}
=\mathbb{F}[I_1,I_2]$.
\end{corollary}
\begin{proof}
The result is a direct consequence of Theorem \ref{theorem}
and \cite[\textsc{Th\'eor\`{e}me} 1-I]{KPV}.
\end{proof}
\section{Normal forms and invariants}
In the series of papers \cite{BK1}, \cite{BK2}, \cite{BK3}, \cite{BK4} and \cite{BK5},
the normal forms for equivalence classes
of trivectors over a $6$-dimensional vector space
over an arbitrary field $\mathbb{F}$ of characteristic
distinct from $2$ and $3$ under the symplectic group,
is given. According to \cite[Theorem 2.1]{BK5}
every non-zero trivector of $V$ is equivalent
with (at least) one of the following trivectors:
\[
\begin{array}
[c]{ll}
\chi _{A_1}=\bar{e}_1^\ast \wedge \bar{e}_2^\ast
\wedge \bar{e}_3^\ast ,
& \chi _{A_2}=\bar{e}_1^\ast \wedge \bar{e}_2^\ast
\wedge \bar{f}_2^\ast ,
\end{array}
\]
\[
\begin{array}
[c]{rl}
\chi _{B_1}= & \bar{e}_1^\ast \wedge \bar{e}_2^\ast
\wedge \bar{e}_3^\ast
+\bar{e}_1^\ast \wedge \bar{f}_1^\ast
\wedge \bar{f}_3^\ast ,\\
\chi _{B_2}= & \bar{e}_1^\ast \wedge \bar{e}_2^\ast
\wedge \bar{f}_2^\ast +\bar{e}_1^\ast \wedge \bar{f}_1^\ast
\wedge \bar{e}_3^\ast,\\
\chi _{B_3}= & \bar{e}_1^\ast \wedge \bar{e}_2^\ast
\wedge \bar{f}_2^\ast +\bar{e}_1^\ast
\wedge \bar{e}_3^\ast \wedge \bar{f}_3^\ast ,\\
\chi _{B_4}(\lambda )= & \bar{e}_1^\ast
\wedge \bar{e}_2^\ast \wedge \bar{e}_3^\ast
+\lambda\bar{e}_1^\ast \wedge \bar{f}_2^\ast
\wedge \bar{f}_3^\ast ,\\
\chi _{B_5}(\lambda )= & \lambda\bar{e}_1^\ast
\wedge \bar{e}_2^\ast \wedge \bar{f}_2^\ast
+\bar{e}_1^\ast \wedge (\bar{e}_2^\ast
-\bar{e}_3^\ast )\wedge (\bar{f}_2^\ast +\bar{f}_3^\ast ),
\end{array}
\]
\[
\begin{array}
[c]{rl}
\chi _{C_1}(\lambda )= & \bar{e}_1^\ast
\wedge \bar{e}_2^\ast \wedge \bar{e}_3^\ast
+\lambda\bar{f}_1^\ast \wedge \bar{f}_2^\ast
\wedge \bar{f}_3^\ast ,\\
\chi _{C_2}(\lambda)= & \bar{f}_1^\ast
\wedge (\bar{e}_2^\ast +\bar{e}_3^\ast )
\wedge (\bar{f}_2^\ast -\bar{f}_3^\ast )
+\lambda \bar{e}_1^\ast \wedge \bar{e}_2^\ast
\wedge \bar{f}_2^\ast ,\\
\chi _{C_3}(\lambda)= & \bar{e}_1^\ast
\wedge \bar{e}_2^\ast \wedge\bar{f}_2^\ast
+\lambda\bar{f}_1^\ast \wedge \bar{e}_3^\ast
\wedge\bar{f}_3^\ast ,\\
\chi _{C_4}(\lambda)= & \bar{f}_1^\ast
\wedge \bar{e}_3^\ast \wedge(\bar{e}_2^\ast
+\bar{f}_3^\ast )+\lambda\bar{e}_1^\ast
\wedge \bar{e}_2^\ast \wedge \bar{f}_2^\ast ,\\
\chi _{C_5}(\lambda )= & \bar{e}_1^\ast
\wedge \bar{e}_3^\ast \wedge
(\bar{f}_2^\ast +\bar{f}_3^\ast )
+\lambda \bar{e}_2^\ast \wedge \bar{f}_3^\ast
\wedge (\bar{f}_1^\ast +\bar{e}_3^\ast ),\\
\chi _{C_6}(\lambda,\varepsilon)= & \bar{f}_1^\ast
\wedge (\bar{e}_2^\ast +\bar{e}_3^\ast )
\wedge (\bar{f}_2^\ast +\varepsilon \bar{f}_3^\ast)
+\lambda \bar{e}_1^\ast \wedge \bar{e}_2^\ast
\wedge \bar{f}_2^\ast ,
\end{array}
\]
\[
\begin{array}
[c]{rl}
\chi _{D_1}= & \bar{e}_1^\ast \wedge \bar{e}_2^\ast
\wedge \bar{f}_2^\ast +\bar{e}_2^\ast
\wedge \bar{f}_1^\ast \wedge \bar{e}_3^\ast
+\bar{f}_1^\ast \wedge \bar{e}_1^\ast
\wedge \bar{f}_3^\ast ,\\
\chi _{D_2}(\lambda)= & \lambda \bar{e}_1^\ast
\wedge \bar{e}_2^\ast \wedge \bar{f}_3^\ast
+\bar{e}_2^\ast \wedge \bar{f}_1^\ast
\wedge \bar{e}_3^\ast +\bar{f}_1^\ast
\wedge \bar{e}_1^\ast \wedge \bar{f}_2^\ast ,\\
\chi _{D_3}(\lambda _1,\lambda _2)= & \bar{e}_1^\ast
\wedge \bar{e}_2^\ast \wedge \bar{f}_3^\ast
+\lambda _1\bar{e}_2^\ast \wedge \bar{e}_3^\ast
\wedge \bar{f}_1^\ast +\lambda _2\bar{e}_3^\ast
\wedge \bar{e}_1^\ast \wedge \bar{f}_2^\ast ,\\
\chi _{D_4}(\lambda _1,\lambda _2)= & \bar{e}_1^\ast
\wedge \bar{e}_2^\ast \wedge \bar{f}_3^\ast
+\lambda _1\bar{e}_2^\ast \wedge \bar{e}_3^\ast
\wedge (\bar{f}_1^\ast +\bar{f}_3^\ast )
+\lambda _2\bar{e}_3^\ast \wedge \bar{e}_1^\ast
\wedge \bar{f}_2^\ast ,\\
\chi _{D_5}(\lambda)= & \bar{e}_1^\ast
\wedge \bar{e}_2^\ast \wedge \bar{f}_3^\ast
+\lambda \bar{e}_2^\ast \wedge \bar{e}_3^\ast
\wedge (\bar{f}_1^\ast +\bar{f}_2^\ast +\bar{f}_3^\ast )
-\bar{e}_3^\ast \wedge \bar{e}_1^\ast
\wedge \bar{f}_2^\ast ,\\
\chi _{D_6}= & -\bar{e}_1^\ast \wedge \bar{e}_2^\ast
\wedge \bar{f}_2^\ast +\bar{e}_2^\ast
\wedge \bar{e}_3^\ast \wedge \bar{f}_1^{\ast
}+\bar{e}_3^\ast \wedge \bar{e}_1^\ast
\wedge \bar{f}_3^\ast ,
\end{array}
\]
\[
\begin{array}
[c]{rl}
\chi _{E_1}(a,b,h_1,h_2,h_3)
= & 2\bar{e}_1^\ast \wedge \bar{e}_2^\ast
\wedge \bar{e}_3^\ast \\
& \multicolumn{1}{r}{+a\left( h_1\bar{f}_1^\ast
\wedge \bar{e}_2^\ast \wedge \bar{e}_3^\ast
+h_2\bar{e}_1^\ast \wedge \bar{f}_2^\ast
\wedge \bar{e}_3^\ast +h_3\bar{e}_1^\ast
\wedge \bar{e}_2^\ast \wedge \bar{f}_3^\ast \right) }\\
& \multicolumn{1}{r}{+(a^2+2b)\left( h_1h_2\bar{f}_1^\ast
\wedge \bar{f}_2^\ast
\wedge \bar{e}_3^\ast +h_1h_3\bar{f}_1^\ast
\wedge \bar{e}_2^\ast \wedge \bar{f}_3^\ast \right. }\\
& \multicolumn{1}{r}{\left.  +h_2h_3\bar{e}_1^\ast
\wedge \bar{f}_2^\ast \wedge \bar{f}_3^\ast \right)
+h_1h_2h_3(a^2+3b)\bar{f}_1^\ast \wedge \bar{f}_2^\ast
\wedge \bar{f}_3^\ast ,}
\end{array}
\]
\[
\begin{array}
[c]{rl}
\chi _{E_2}(a,b,k)= & \bar{e}_1^\ast \wedge \bar{e}_2^\ast
\wedge \bar{f}_2^\ast +\bar{e}_1^\ast \wedge \bar{e}_3^\ast
\wedge \bar{f}_3^\ast \\
& \multicolumn{1}{r}{+k\left( \bar{f}_1^\ast
\wedge \bar{e}_2^\ast \wedge \bar{f}_3^\ast
-b\bar{f}_1^\ast \wedge \bar{f}_2^\ast
\wedge \bar{e}_3^\ast +a\bar{f}_1^\ast
\wedge \bar{e}_3^\ast \wedge
\bar{f}_3^\ast \right) ,\medskip}\\
\chi _{E_3}(a,b,k,h)= & \bar{e}_1^\ast
\wedge \bar{e}_2^\ast \wedge \bar{f}_2^\ast
+\bar{e}_1^\ast \wedge \bar{e}_3^\ast
\wedge \bar{f}_3^\ast +k\left( \bar{f}_1^\ast
\wedge \bar{e}_2^\ast \wedge \bar{f}_3^\ast
-b\bar{f}_1^\ast \wedge \bar{f}_2^\ast
\wedge \bar{e}_3^\ast \right. \\
& \multicolumn{1}{r}{\left. +a\bar{f}_1^\ast
\wedge \bar{e}_3^\ast \wedge \bar{f}_3^\ast
\right) +h\bar{e}_1^\ast \wedge \bar{f}_2^\ast
\wedge \bar{f}_3^\ast ,}
\end{array}
\]
\[
\begin{array}
[c]{rl}
\chi _{E_4}(a,b,k,h_1,h_2)= & \left[  1-h_1h_2(a^2+4b)\right]
\bar{e}_1^\ast \wedge \bar{e}_2^\ast \wedge \bar{f}_2^\ast \\
& +\left[  1+h_1h_2(a^2+4b)\right] \bar{e}_1^\ast
\wedge \bar{e}_3^\ast \wedge \bar{f}_3^\ast \\
& +k\left[  \bar{f}_1^\ast \wedge \bar{e}_2^\ast
\wedge \bar{f}_3^\ast -b(1-h_1h_2(a^2+4b))
\bar{f}_1^\ast \wedge \bar{f}_2^\ast
\wedge \bar{e}_3^\ast \right. \\
& \left.  +a\bar{f}_1^\ast \wedge \bar{e}_3^\ast
\wedge \bar{f}_3^\ast \right]
+h_1\left[ 1-h_1h_2(a^2+4b)\right]  \bar{e}_1^\ast
\wedge \bar{f}_2^\ast \wedge \bar{f}_3^\ast \\
& \multicolumn{1}{r}{+(a^2+4b)h_2\bar{e}_1^\ast
\wedge \bar{e}_2^\ast \wedge \bar{e}_3^\ast ,}\\
& h_1h_2(a^2+4b)\neq 1,
\end{array}
\]
\[
\begin{array}
[c]{rl}
\chi _{E_5}(a,b,k)= & \bar{f}_1^\ast \wedge \bar{e}_2^\ast
\wedge \bar{f}_3^\ast +2\bar{e}_1^\ast
\wedge \bar{f}_1^\ast \wedge \bar{e}_2^\ast
-a\bar{f}_1^\ast \wedge \bar{e}_2^\ast
\wedge \bar{f}_2^\ast \\
& \multicolumn{1}{r}{+a\bar{f}_1^\ast
\wedge \bar{e}_3^\ast \wedge \bar{f}_3^\ast
+a\bar{e}_1^\ast \wedge \bar{f}_1^\ast
\wedge \bar{e}_3^\ast +(a^2+b)\bar{f}_1^\ast
\wedge \bar{f}_2^\ast \wedge \bar{e}_3^\ast }\\
& \multicolumn{1}{r}{+k\left( a\bar{e}_1^\ast
\wedge \bar{f}_2^{\ast
}\wedge \bar{e}_3^\ast -\bar{e}_1^\ast
\wedge \bar{e}_2^\ast \wedge \bar{f}_2^\ast
+\bar{e}_1^\ast \wedge \bar{e}_3^\ast
\wedge \bar{f}_3^\ast \right) ,}
\end{array}
\]
(As we are assuming that $\mathbb{F}$ has characteristic
distinct from $2$ and $3$, we do not consider the trivectors
of types $(E1^\prime )$, $(E2^\prime )$, and $(E3^\prime )$,
because they only exist when the characteristic is $2$; indeed,
for such trivectors one needs separable quadratic extensions
of $\mathbb{F}$, which can only exist in the case
of characteristic $2$.)

With the same notations as above, we have
$\bar{e}_1^\ast =v^1$,
$\bar{e}_2^\ast =v^2$, $\bar{e}_3^\ast =v^3$,
$\bar{f}_1^\ast =v^4$, $\bar{f}_2^\ast =v^{5}$,
$\bar{f}_3^\ast =v^6$.

In this section we study the evaluation map
$(I_1,I_2)\colon \wedge ^3V^\ast /GL(V)
\to \mathbb{F}^2$,
by using the normal forms given in \cite{BK5}.
This throws light on the quotient space. We have
\[
\begin{array}{llll}
I_j(\chi _{A_h})=0, \quad I_j(\chi _{B_i})=0,
& h,j=1,2, & 1\leq i\leq 5,\\
I_1(\chi _{C_h})=0, & I_2(\chi _{C_h})=-2^3\cdot 3^2\lambda ^2,
& 1\leq h\leq 2, &  \\
I_1(\chi _{C_h})=\lambda ^2, & I_2(\chi _{C_h})
=2^4\cdot 3\lambda ^2, & 3\leq h\leq 4, &  \\
I_1(\chi _{C_5})=0, & I_2(\chi _{C_5})
=-2^3\cdot 3^2\lambda ^2, &  &  \\
I_1(\chi _{C_6})=\lambda ^2
\varepsilon (\varepsilon +1),
& I_2(\chi _{C_6})=2^3\cdot 3\lambda ^2
\varepsilon (2\varepsilon +5), &
\varepsilon \neq 0,-1, &  \\
I_j(\chi _{D_h})=0, & h\in \{ 1,3,4,5,6\} ,
& 1\leq j\leq 2, &  \\
I_1(\chi _{D_2})=\lambda , & I_2(\chi _{D_2})
=2^3\cdot 3\cdot 5\lambda , &  &  \\
I_1(\chi _{E_1})=0, &  &  &
\end{array}
\]
\begin{align*}
I_2(\chi _{E_1})&=-2^3\cdot 3^2(h_1h_2h_3)^2
\left( 2^2 (a^2+3b)^2(1-2a)+(a^2+2b)^2(5a^2+2^4 b)
\right) , \quad \\
I_1(\chi _{E_{j}})&=k^2(a^2+4 b),
\qquad \qquad I_2(\chi _{E_j})
=2^4\cdot 3k^2(a^2+4b), \quad 2\leq j\leq 3,\\
I_1(\chi _{E_4}) &
=k^2 (a^2+4b)\left( 1-(a^2+4 b) h_1 h_2 \right) ,
\qquad \qquad \qquad \qquad\\
I_2(\chi _{E_4}) &
=2^3\cdot 3 k^2 (a^2+4 b) \left( 1-(a^2+4 b)h_1 h_2 \right)
(2+3 (a^2 + 4 b) h_1 h_2), \\
I_1(\chi _{E_5})&=0,
\qquad I_2(\chi _{E_5})=-2^3 \cdot 3^2 k^2 (a^2+4 b)
\qquad \qquad \qquad
\end{align*}
When $\mathbb{F}$ is an algebraically closed field
of characteristic distinct from $2$, Popov \cite{P}
gave another alternative for the normal forms for
equivalence classes of trivectors over
a $6$-dimensional vector space over $\mathbb{F}$
under the symplectic group. We list these normal forms,
following the notations in \cite[Theorem 2.1]{BK5},
with $q,p\in\mathbb{F}^\ast $:
\[
\begin{array}
[c]{rl}
\chi _{P_1}= & 0,\\
\chi _{P_2}= & \bar{e}_1^\ast \wedge \bar{e}_2^\ast
\wedge \bar{f}_2^\ast +\bar{e}_1^\ast
\wedge \bar{e}_3^\ast \wedge \bar{f}_3^\ast ,\\
\chi _{P_3}(q)= & \bar{e}_1^\ast \wedge
\bar{e}_2^\ast \wedge \bar{e}_3^\ast
+q\bar{f}_1^\ast \wedge \bar{f}_2^\ast
\wedge \bar{f}_3^\ast ,\\
\chi _{P_4}(q)= & \bar{e}_1^\ast
\wedge \bar{e}_2^\ast \wedge \bar{e}_3^\ast
+q\bar{f}_1^\ast \wedge \bar{f}_2^\ast
\wedge \bar{f}_3^\ast +\bar{e}_1^\ast
\wedge \bar{e}_2^\ast \wedge \bar{f}_2^\ast
+\bar{e}_1^\ast \wedge \bar{e}_3^\ast
\wedge \bar{f}_3^\ast ,\\
\chi _{P_5}(q)= & \bar{e}_1^\ast
\wedge \bar{e}_2^\ast \wedge \bar{e}
_3^\ast +q\bar{f}_1^\ast \wedge
\bar{f}_2^\ast \wedge \bar{f}_3^\ast
+\bar{e}_1^\ast \wedge \bar{e}_2^\ast
\wedge \bar{f}_2^\ast +\bar{e}_1^\ast
\wedge \bar{e}_3^\ast \wedge \bar{f}_3^\ast \\
& +\bar{f}_2^\ast \wedge \bar{e}_1^\ast
\wedge \bar{f}_1^\ast +\bar{f}_2^\ast
\wedge \bar{e}_3^\ast \wedge \bar{f}_3^\ast ,\\
\chi _{P_6}(q,p)= & \bar{e}_1^\ast
\wedge \bar{e}_2^\ast \wedge \bar{e}_3^\ast
+q\bar{f}_1^\ast \wedge \bar{f}_2^\ast
\wedge \bar{f}_3^\ast +\bar{e}_1^\ast
\wedge \bar{e}_2^\ast \wedge \bar{f}_2^\ast
+\bar{e}_1^\ast \wedge \bar{e}_3^\ast
\wedge \bar{f}_3^\ast \\
& +p\bar{f}_1^\ast \wedge \bar{e}_2^\ast
\wedge \bar{f}_2^\ast
+p\bar{f}_1^\ast \wedge \bar{e}_3^\ast
\wedge \bar{f}_3^\ast ,\\
\chi _{P_{7}}= & \bar{e}_1^\ast
\wedge \bar{e}_2^\ast \wedge \bar{e}_3^\ast ,\\
\chi _{P_{8}}= & \bar{e}_1^\ast
\wedge \bar{e}_2^\ast \wedge \bar{e}_3^\ast
+\bar{e}_1^\ast \wedge \bar{e}_2^\ast
\wedge \bar{f}_2^\ast +\bar{e}_1^\ast
\wedge \bar{e}_3^\ast \wedge \bar{f}_3^\ast ,\\
\chi _{P_{9}}= & \bar{e}_1^\ast
\wedge \bar{e}_2^\ast \wedge \bar{e}
_3^\ast +\bar{f}_1^\ast \wedge \bar{e}_2^\ast
\wedge \bar{f}_2^\ast+\bar{f}_1^\ast
\wedge \bar{e}_3^\ast \wedge \bar{f}_3^\ast ,
\end{array}
\]
\[
\begin{array}
[c]{rl}
\chi _{P_{10}}= & \bar{e}_1^\ast \wedge \bar{f}_1^\ast
\wedge \bar{e}_3^\ast +\bar{f}_2^\ast
\wedge \bar{e}_2^\ast \wedge \bar{e}_3^\ast
+\bar{e}_1^\ast \wedge \bar{e}_2^\ast
\wedge \bar{f}_3^\ast ,\\
\chi _{P_{11}}= & \bar{e}_1^\ast \wedge \bar{f}_1^\ast
\wedge \bar{e}_3^\ast +\bar{f}_2^\ast
\wedge \bar{e}_2^\ast \wedge \bar{e}_3^\ast
+\bar{e}_1^\ast \wedge \bar{e}_2^\ast
\wedge \bar{f}_3^\ast +\bar{e}_1^\ast
\wedge \bar{e}_2^\ast \wedge \bar{f}_2^\ast \\
& +\bar{e}_1^\ast \wedge \bar{e}_3^\ast
\wedge \bar{f}_3^\ast ,\\
\chi _{P_{12}}(q)= & \bar{e}_1^\ast
\wedge \bar{f}_1^\ast \wedge \bar{e}_3^\ast
+\bar{f}_2^\ast \wedge \bar{e}_2^\ast \wedge \bar{e}
_3^\ast +\bar{e}_1^\ast \wedge \bar{e}_2^\ast
\wedge \bar{f}_3^\ast +q\bar{e}_3^\ast
\wedge \bar{e}_1^\ast \wedge \bar{f}_1^\ast \\
& \multicolumn{1}{r}{+q\bar{e}_3^\ast
\wedge \bar{e}_2^\ast \wedge \bar{f}_2^\ast ,}\\
\chi _{P_{13}}= & \bar{e}_1^\ast
\wedge \bar{f}_1^\ast \wedge \bar{e}_3^\ast
+\bar{f}_2^\ast \wedge \bar{e}_2^\ast
\wedge \bar{e}_3^\ast +\bar{e}_1^\ast
\wedge \bar{e}_2^\ast \wedge \bar{f}_3^\ast
+\bar{f}_1^\ast \wedge \bar{e}_2^\ast
\wedge \bar{f}_2^\ast \\
& \multicolumn{1}{r}{+\bar{f}_1^\ast
\wedge \bar{e}_3^\ast \wedge \bar{f}_3^\ast ,}\\
\chi _{P_{14}}(q)= & \bar{e}_1^\ast
\wedge \bar{f}_1^\ast \wedge \bar
{e}_3^\ast +\bar{f}_2^\ast
\wedge \bar{e}_2^\ast \wedge \bar{e}_3^\ast
+\bar{e}_1^\ast \wedge \bar{e}_2^\ast
\wedge \bar{f}_3^\ast +q\bar{f}_3^\ast
\wedge \bar{e}_1^\ast \wedge \bar{f}_1^\ast \\
& \multicolumn{1}{r}{+q\bar{f}_3^\ast
\wedge \bar{e}_2^\ast \wedge \bar{f}_2^\ast ,}
\end{array}
\]
\[
\begin{array}
[c]{rl}
\chi _{P_{15}}= & \bar{e}_1^\ast
\wedge \bar{f}_1^\ast \wedge \bar{e}_3^\ast
+\bar{f}_2^\ast \wedge \bar{e}_2^\ast
\wedge \bar{e}_3^\ast ,\\
\chi _{P_{16}}(q)= & \bar{e}_1^\ast
\wedge \bar{f}_1^\ast \wedge \bar{e}_3^\ast
+\bar{f}_2^\ast \wedge \bar{e}_2^\ast
\wedge \bar{e}_3^\ast +q\bar{f}_3^\ast
\wedge \bar{e}_1^\ast \wedge \bar{f}_1^\ast
+q\bar{f}_3^\ast \wedge \bar{e}_2^\ast
\wedge \bar{f}_2^\ast ,\\
\chi _{P_{17}}(q)= & \bar{e}_1^\ast
\wedge \bar{f}_1^\ast \wedge \bar{e}_3^\ast
+\bar{f}_2^\ast \wedge \bar{e}_2^\ast
\wedge \bar{e}_3^\ast +q\bar{e}_3^\ast
\wedge \bar{e}_1^\ast \wedge \bar{f}_1^\ast
+q\bar{e}_3^\ast \wedge \bar{e}_2^\ast
\wedge \bar{f}_2^\ast ,\\
\chi _{P_{18}}= & \bar{e}_1^\ast
\wedge \bar{f}_1^\ast \wedge \bar{e}_3^\ast
+\bar{f}_2^\ast \wedge \bar{e}_2^\ast
\wedge \bar{e}_3^\ast
+\bar{e}_1^\ast \wedge \bar{e}_2^\ast
\wedge \bar{f}_2^\ast +\bar{e}_1^\ast
\wedge \bar{e}_2^\ast \wedge \bar{f}_3^\ast ,\\
\chi _{P_{19}}= & \bar{e}_1^\ast
\wedge \bar{f}_1^\ast \wedge \bar{e}_3^\ast
+\bar{f}_2^\ast \wedge \bar{e}_2^\ast
\wedge \bar{e}_3^{\ast
}+\bar{e}_1^\ast \wedge \bar{e}_2^\ast
\wedge \bar{f}_2^\ast +\bar{e}_1^\ast
\wedge \bar{e}_2^\ast \wedge \bar{f}_3^\ast \\
& +\bar{e}_2^\ast \wedge \bar{e}_1^\ast
\wedge \bar{f}_1^\ast +\bar{e}_2^\ast
\wedge \bar{e}_3^\ast \wedge \bar{f}_3^\ast .
\end{array}
\]

Then, we obtain
\[
\begin{array}
[c]{lll}
I_{j}(\chi _{P_h})=0, & 1\leq j\leq2,
& h\in \{ 2,7,8,9,10,11,12,13,15,17,18,19\} ,
\end{array}
\]
\[
\begin{array}
[c]{lll}
I_1(\chi _{P_h})=0, & I_2(\chi _{P_h})=-2^3\cdot 3^2q^2,
& h\in\{ 3,4,5\} ,\\
I_1(\chi _{P_6})=-2^2pq, & I_2(\chi _{P_6})
=-2^3\cdot 3q(3q+2^3p), & \\
I_1(\chi _{P_h})=2^2q^2,
& I_2(\chi _{P_h})=2^6\cdot 3 q^2,
& h\in \{ 14,16\} .
\end{array}
\]

\subsubsection*{Acknowledgment}
The authors wish to thank the unknown referee for his careful and thorough read of our manuscript and his valuable comments.

\section{Appendix}
\begin{align*}
I_1  & =y_{135}y_{234}y_{256}^2-y_{126}y_{234}y_{356}^2-y_{134}
y_{136}y_{236}y_{456}-y_{126}y_{134}y_{145}y_{346}\\
& +2y_{126}y_{135}y_{245}y_{346}+y_{124}y_{135}y_{245}y_{256}-y_{126}
y_{145}y_{235}y_{245}\\
& -3y_{126}y_{145}y_{235}y_{346}+y_{124}y_{135}y_{146}y_{346}+3y_{124}
y_{135}y_{256}y_{346}\\
& -y_{124}y_{135}y_{146}y_{245}+y_{145}y_{156}y_{234}y_{236}+2y_{135}
y_{146}y_{234}y_{256}\\
& +y_{124}y_{156}y_{234}y_{356}+y_{156}y_{235}y_{236}y_{346}-y_{156}
y_{235}y_{236}y_{245}\\
& -y_{146}y_{156}y_{234}y_{235}+2y_{126}y_{145}y_{234}y_{356}-y_{156}
y_{234}y_{236}y_{356}\\
& +y_{136}y_{156}y_{234}y_{346}-y_{156}y_{234}y_{235}y_{256}+y_{136}
y_{156}y_{234}y_{245}\\
& -y_{124}y_{156}y_{235}y_{245}+y_{146}^2y_{235}^2+y_{145}^2y_{236}
^2+y_{124}^2y_{356}^2\\
& +2y_{124}y_{156}y_{236}y_{345}-3y_{124}y_{156}y_{235}y_{346}-y_{125}
y_{126}y_{345}y_{346}\\
& +y_{123}y_{136}y_{346}y_{456}-y_{123}y_{145}y_{236}y_{456}-y_{125}
y_{126}y_{245}y_{345}\\
& +2y_{125}y_{146}y_{235}y_{346}-2y_{136}y_{146}y_{235}y_{245}-2y_{136}
y_{145}y_{236}y_{245}\\
& -3y_{123}y_{146}y_{245}y_{356}-2y_{145}y_{146}y_{235}y_{236}-3y_{125}
y_{134}y_{236}y_{456}\\
& +2y_{125}y_{136}y_{234}y_{456}+2y_{123}y_{156}y_{245}y_{346}+2y_{123}
y_{145}y_{246}y_{356}\\
& -3y_{125}y_{146}y_{236}y_{345}-y_{123}y_{145}y_{146}y_{346}+y_{123}
y_{145}y_{146}y_{245}\\
& -y_{123}y_{145}y_{245}y_{256}-y_{123}y_{146}y_{346}y_{356}-y_{125}
y_{234}y_{256}y_{356}\\
& +y_{134}y_{136}y_{145}y_{246}+y_{134}y_{136}y_{246}y_{356}-y_{125}
y_{145}y_{234}y_{256}\\
& -3y_{125}y_{146}y_{234}y_{356}-3y_{123}y_{145}y_{256}y_{346}+y_{125}
y_{145}y_{146}y_{234}\\
& +y_{123}y_{256}y_{346}y_{356}-y_{123}y_{245}y_{256}y_{356}+y_{123}
y_{136}y_{245}y_{456}\\
& +y_{123}y_{134}y_{256}y_{456}+y_{135}y_{236}y_{245}y_{256}+2y_{134}
y_{136}y_{245}y_{256}\\
& +y_{136}^2y_{245}^2+y_{125}^2y_{346}^2-y_{125}y_{236}y_{256}
y_{345}-y_{135}y_{236}y_{256}y_{346}\\
& +3y_{135}y_{146}y_{236}y_{245}+y_{123}y_{235}y_{256}y_{456}-y_{125}
y_{235}y_{236}y_{456}\\
& +y_{135}y_{146}y_{236}y_{346}+y_{123}y_{236}y_{356}y_{456}-3y_{124}
y_{136}y_{235}y_{456}\\
& +y_{124}y_{125}y_{134}y_{456}+2y_{124}y_{136}y_{245}y_{356}-2y_{124}
y_{145}y_{236}y_{356}\\
& -2y_{124}y_{134}y_{256}y_{356}+y_{125}y_{235}y_{246}y_{356}-y_{125}
y_{134}y_{145}y_{246}\\
& +3y_{125}y_{134}y_{246}y_{356}+2y_{124}y_{135}y_{236}y_{456}-3y_{124}
y_{136}y_{256}y_{345}\\
& +y_{125}y_{145}y_{235}y_{246}+y_{123}y_{146}y_{235}y_{456}+2y_{123}
y_{146}y_{256}y_{345}\\
& +y_{123}y_{125}y_{346}y_{456}-y_{124}y_{136}y_{146}y_{345}+y_{123}
y_{134}y_{146}y_{456}\\
& +y_{136}y_{236}y_{256}y_{345}+y_{136}y_{235}y_{236}y_{456}-y_{124}
y_{126}y_{145}y_{345}\\
& +y_{125}y_{156}y_{234}y_{245}-y_{125}y_{135}y_{246}y_{346}+y_{126}
y_{235}y_{346}y_{356}\\
& -y_{125}y_{135}y_{245}y_{246}-2y_{125}y_{136}y_{245}y_{346}-2y_{134}
y_{146}y_{235}y_{256}
\end{align*}
\begin{align*}
& +y_{134}^2y_{256}^2+2y_{134}y_{145}y_{236}y_{256}-y_{135}y_{235}
y_{246}y_{256}+y_{136}y_{234}y_{256}y_{356}\\
& -y_{136}y_{145}y_{146}y_{234}-y_{126}y_{235}y_{245}y_{356}-y_{136}
y_{146}y_{236}y_{345}\\
& -y_{136}y_{146}y_{234}y_{356}-3y_{136}y_{145}y_{234}y_{256}+y_{125}
y_{156}y_{234}y_{346}\\
& +y_{126}y_{145}y_{236}y_{345}-y_{126}y_{136}y_{245}y_{345}+y_{126}
y_{134}y_{146}y_{345}
\end{align*}
\begin{align*}
& +2y_{125}y_{136}y_{246}y_{345}+y_{126}y_{134}y_{256}y_{345}+y_{126}
y_{235}y_{256}y_{345}\\
& -y_{134}y_{135}y_{246}y_{256}-y_{126}y_{236}y_{345}y_{356}-y_{124}
y_{135}y_{246}y_{356}\\
& -y_{124}y_{145}y_{156}y_{234}+y_{126}y_{146}y_{235}y_{345}-y_{126}
y_{136}y_{345}y_{346}\\
& +y_{124}y_{135}y_{145}y_{246}-y_{124}y_{134}y_{156}y_{346}-y_{134}
y_{156}y_{234}y_{256}\\
& +2y_{134}y_{156}y_{235}y_{246}+2y_{125}y_{145}y_{236}y_{346}-2y_{125}
y_{134}y_{256}y_{346}\\
& +2y_{124}y_{146}y_{235}y_{356}-3y_{134}y_{156}y_{236}y_{245}-y_{134}
y_{135}y_{146}y_{246}\\
& -y_{135}y_{136}y_{245}y_{246}-y_{134}y_{146}y_{156}y_{234}-y_{135}
y_{145}y_{236}y_{246}\\
& -y_{135}y_{136}y_{246}y_{346}-y_{134}y_{156}y_{236}y_{346}+2y_{126}
y_{134}y_{235}y_{456}\\
& -2y_{124}y_{125}y_{346}y_{356}+3y_{136}y_{145}y_{235}y_{246}+y_{126}
y_{134}y_{145}y_{245}\\
& -3y_{126}y_{134}y_{245}y_{356}-y_{126}y_{134}y_{346}y_{356}-y_{124}
y_{134}y_{136}y_{456}\\
& -y_{124}y_{125}y_{235}y_{456}+y_{124}y_{125}y_{146}y_{345}-y_{124}
y_{125}y_{256}y_{345}\\
& -y_{136}y_{235}y_{246}y_{356}+y_{123}y_{124}y_{145}y_{456}-y_{123}
y_{124}y_{356}y_{456}\\
& +y_{125}^2y_{246}y_{345}+y_{123}y_{256}^2y_{345}+y_{123}y_{156}
y_{245}^2+y_{126}y_{135}y_{245}^2\\
& -y_{135}y_{236}^2y_{456}+y_{135}y_{146}^2y_{234}+y_{125}^2
y_{234}y_{456}-y_{123}y_{145}^2y_{246}\\
& -y_{156}y_{236}^2y_{345}+y_{156}y_{235}^2y_{246}+y_{134}^2
y_{156}y_{246}+y_{123}y_{156}y_{346}^2\\
& +y_{126}y_{135}y_{346}^2-y_{124}^2y_{156}y_{345}+y_{123}y_{146}
^2y_{345}+y_{126}y_{235}^2y_{456}\\
& +y_{126}y_{134}^2y_{456}-y_{124}^2y_{135}y_{456}+y_{136}^2
y_{246}y_{345}-y_{123}y_{246}y_{356}^2\\
& -y_{126}y_{145}^2y_{234}+y_{136}^2y_{234}y_{456}-y_{135}y_{146}
y_{235}y_{246}\\
& +y_{123}y_{125}y_{245}y_{456}+y_{124}y_{134}y_{156}y_{245}+y_{135}
y_{236}y_{246}y_{356}\\
& +y_{124}y_{126}y_{345}y_{356}.
\end{align*}
\begin{align*}
I_2 &  =24\left( -5y_{124}y_{136}y_{146}y_{345}-5y_{126}y_{145}
y_{235}y_{245}-5y_{124}y_{156}y_{235}y_{245}\right. \\
&  -4y_{136}y_{146}y_{235}y_{245}-5y_{136}y_{146}y_{236}y_{345}-5y_{136}
y_{146}y_{234}y_{356}\\
&  -5y_{124}y_{125}y_{256}y_{345}-5y_{125}y_{126}y_{245}y_{345}-2y_{125}
y_{136}y_{246}y_{345}\\
&  +y_{125}y_{126}y_{345}y_{346}-5y_{124}y_{135}y_{146}y_{245}+5y_{124}
y_{135}y_{245}y_{256}\\
&  -5y_{125}y_{135}y_{245}y_{246}+y_{135}y_{136}y_{245}y_{246}-2y_{126}
y_{135}y_{245}y_{346}\\
&  +4y_{124}y_{146}y_{235}y_{356}-5y_{126}y_{235}y_{245}y_{356}-5y_{136}
y_{235}y_{246}y_{356}\\
&  +6y_{123}y_{156}y_{234}y_{456}+6y_{135}y_{156}y_{234}y_{246}+6y_{126}
y_{156}y_{234}y_{345}\\
&  +5y_{126}y_{134}y_{146}y_{345}+6y_{123}y_{126}y_{345}y_{456}+6y_{126}
y_{135}y_{246}y_{345}\\
&  -y_{123}y_{136}y_{245}y_{456}+4y_{124}y_{136}y_{245}y_{356}+y_{126}
y_{136}y_{245}y_{345}\\
&  -4y_{134}y_{146}y_{235}y_{256}+12y_{123}y_{156}y_{246}y_{345}
+12y_{126}y_{135}y_{234}y_{456}\\
&  +5y_{123}y_{146}y_{145}y_{245}-5y_{123}y_{145}y_{245}y_{256}-2y_{126}
y_{134}y_{235}y_{456}\\
&  +5y_{136}y_{235y_{236}}y_{456}-3y_{124}y_{136}y_{235}y_{456}-5y_{125}
y_{235}y_{236}y_{456}\\
&  -2y_{123}y_{146}y_{256}y_{345}-3y_{123}y_{146}y_{245}y_{356}-5y_{123}
y_{245}y_{256}y_{356}\\
&  -5y_{134}y_{136}y_{236}y_{456}+5y_{123}y_{256}y_{346}y_{356}-5y_{124}
y_{134}y_{136}y_{456}\\
&  -2y_{124}y_{135}y_{236}y_{456}-2y_{123}y_{145}y_{246}y_{356}+5y_{124}
y_{125}y_{134}y_{456}\\
&  -3y_{123}y_{145}y_{256}y_{346}-3y_{125}y_{134}y_{236}y_{456}-2y_{123}
y_{156}y_{245}y_{346}\\
&  -2y_{125}y_{136}y_{234}y_{456}-5y_{124}y_{125}y_{456}y_{235}+5y_{123}
y_{125}y_{456}y_{245}\\
&  -4y_{124}y_{134}y_{256}y_{356}-5y_{123}y_{145}y_{146}y_{346}-5y_{123}
y_{146}y_{346}y_{356}\\
&  +5y_{123}y_{136}y_{456}y_{346}-5y_{126}(y_{145})^2y_{234}-5y_{156}
(y_{236})^2y_{345}\\
&  -5(y_{124})^2y_{156}y_{345}-5y_{126}y_{234}(y_{356})^2+5y_{135}
(y_{146})^2y_{234}\\
&  +\!5y_{156}(y_{235})^2y_{246}+5(y_{134})^2y_{156}y_{246}+2(y_{124}
)^2(y_{356})^2\\
&  -4(y_{156})^2(y_{234})^2-3(y_{126})^2(y_{345})^2-3(y_{123}
)^2(y_{456})^2\\
&  -3(y_{135})^2(y_{246})^2+2(y_{136})^2(y_{245})^2+5y_{135}
y_{234}(y_{256})^2\\
&  +5y_{125}^2y_{345}y_{246}+5y_{135}y_{245}^2y_{126}+5y_{126}y_{346}
^2y_{135}\\
&  +5(y_{136})^2y_{246}y_{345}+5y_{123}y_{156}(y_{245})^2+5(y_{136}
)^2y_{234}y_{456}\\
&  +5y_{126}(y_{235})^2y_{456}+5y_{123}(y_{256})^2y_{345}-5y_{123}
y_{246}(y_{356})^2\\
&  -5y_{135}(y_{236})^2y_{456}-5y_{123}(y_{145})^2y_{246}-5(y_{124}
)^2y_{135}y_{456}\\
&  +5(y_{125})^2y_{234}y_{456}+5y_{126}(y_{134})^2y_{456}+5y_{123}
(y_{146})^2y_{345}\\
&  +5y_{123}y_{156}(y_{346})^2+5y_{123}y_{235}y_{256}y_{456}+5y_{126}
y_{256}y_{235}y_{345}\\
&  +5y_{123}y_{134}y_{146}y_{456}+y_{123}y_{145}y_{236}y_{456}+y_{123}
y_{124}y_{356}y_{456}
\end{align*}
\begin{align*}
&  +6y_{123}y_{135}y_{246}y_{456}+y_{135}y_{145}y_{236}y_{246}+y_{124}
y_{135}y_{246}y_{356}\\
&  +4y_{125}y_{145}y_{236}y_{346}-y_{123}y_{125}y_{346}y_{456}+y_{125}
y_{135}y_{246}y_{346}\\
&  -5y_{125}y_{134}y_{145}y_{246}+5y_{125}y_{145}y_{235}y_{246}-5y_{125}
y_{145}y_{234}y_{256}\\
&  -5y_{134}y_{156}y_{236}y_{346}+5y_{135}y_{146}y_{236}y_{346}-4y_{125}
y_{134}y_{256}y_{346}\\
&  +5y_{124}y_{135}y_{145}y_{246}-4y_{124}y_{145}y_{236}y_{356}+5y_{123}
y_{236}y_{356}y_{456}\\
&  +5y_{135}y_{236}y_{246}y_{356}+5y_{123}y_{124}y_{145}y_{456}-y_{123}
y_{146}y_{235}y_{456}\\
&  -y_{123}y_{134}y_{256}y_{456}-y_{126}y_{146}y_{235}y_{345}-y_{126}
y_{134}y_{256}y_{345}\\
&  -5y_{126}y_{134}y_{145}y_{346}-5y_{126}y_{136}y_{345}y_{346}-5y_{126}
y_{134}y_{346}y_{356}\\
&  -5y_{135}y_{136}y_{246}y_{346}+5y_{134}y_{136}y_{246}y_{356}+4y_{134}
y_{145}y_{236}y_{256}\\
&  -5y_{125}y_{236}y_{256}y_{345}-5y_{135}y_{236}y_{256}y_{346}+5y_{125}
y_{156}y_{234}y_{245}\\
&  -y_{136}y_{156}y_{234}y_{245}-y_{125}y_{156}y_{234}y_{346}+5y_{136}
y_{156}y_{234}y_{346}\\
&  +5y_{125}y_{145}y_{146}y_{234}-y_{145}y_{156}y_{234}y_{236}-5y_{124}
y_{145}y_{156}y_{234}\\
&  -3y_{136}y_{145}y_{234}y_{256}-2y_{126}y_{145}y_{234}y_{356}-3y_{125}
y_{146}y_{236}y_{345}\\
&  -y_{126}y_{145}y_{236}y_{345}-2y_{124}y_{156}y_{236}y_{345}+5y_{136}
y_{236}y_{345}y_{256}\\
&  -5y_{126}y_{236}y_{345}y_{356}+4y_{125}y_{146}y_{235}y_{346}-3y_{126}
y_{145}y_{235}y_{346}\\
&  +5y_{156}y_{236}y_{235}y_{346}-3y_{124}y_{156}y_{235}y_{346}+5y_{126}
y_{235}y_{346}y_{356}\\
&  +5y_{245}y_{134}y_{126}y_{145}-3y_{245}y_{134}y_{236}y_{156}+5y_{245}
y_{134}y_{156}y_{124}\\
&  +4y_{245}y_{134}y_{136}y_{256}-3y_{245}y_{134}y_{126}y_{356}+5y_{345}
y_{124}y_{146}y_{125}\\
&  -5y_{345}y_{124}y_{126}y_{145}-3y_{345}y_{124}y_{136}y_{256}-y_{345}
y_{124}y_{126}y_{356}\\
&  -3y_{356}y_{234}y_{146}y_{125}-5y_{356}y_{234}y_{236}y_{156}-y_{356}
y_{234}y_{156}y_{124}\\
&  +5y_{356}y_{234}y_{136}y_{256}-5y_{146}y_{234}y_{136}y_{145}+y_{146}
y_{234}y_{156}y_{235}\\
&  -5y_{146}y_{234}y_{156}y_{134}-2y_{146}y_{234}y_{256}y_{135}+y_{235}
y_{246}y_{135}y_{146}\\
&  +3y_{235}y_{246}y_{136}y_{145}-2y_{235}y_{246}y_{156}y_{134}-5y_{235}
y_{246}y_{256}y_{135}\\
&  +5y_{125}y_{235}y_{246}y_{356}+3y_{135}y_{146}y_{236}y_{245}-4y_{136}
y_{145}y_{236}y_{245}\\
&  -5y_{156}y_{235}y_{236}y_{245}+5y_{135}y_{236}y_{245}y_{256}-5y_{134}
y_{135}y_{146}y_{246}\\
&  +5y_{134}y_{136}y_{145}y_{246}+y_{134}y_{135}y_{246}y_{256}+3y_{125}
y_{134}y_{246}y_{356}\\
&  +5y_{124}y_{135}y_{146}y_{346}-5y_{124}y_{134}y_{156}y_{346}+3y_{124}
y_{135}y_{346}y_{256}\\
&  -4y_{346}y_{124}y_{356}y_{125}-5y_{256}y_{234}y_{156}y_{235}+y_{134}
y_{156}y_{234}y_{256}\\
&  -5y_{125}y_{234}y_{256}y_{356}-4y_{145}y_{146}y_{236}y_{235}-4y_{125}
y_{136}y_{245}y_{346}\\
&  \left.  +2(y_{146})^2(y_{235})^2+2(y_{145})^2(y_{236})^2
+2(y_{134})^2(y_{256}^2)+2(y_{125})^2(y_{346})^2\right) .
\end{align*}

\noindent\textbf{Authors' addresses}

\smallskip

\noindent(J.M.M.) \textsc{Instituto de Tecnolog\'{\i}as
F\'{\i}sicas y de la Informaci\'on, CSIC, C/ Serrano 144,
28006-Madrid, Spain.}

\noindent\emph{E-mail:\/} \texttt{jaime@iec.csic.es}

\medskip

\noindent(L.M.P.C.) \textsc{Departamento de Matem\'atica
Aplicada a las T.I.C., E.T.S.I. Sistemas Inform\'aticos,
Universidad Polit\'ecnica de Madrid,
Carrete\-ra de Valencia, Km.\ 7, 28031-Madrid, Spain.}

\noindent\emph{E-mail:\/} \texttt{luispozo@etsisi.upm.es}

\begin{thebibliography}{9}
\bibitem {BK1} B. De Bruyn, M. Kwiatkowski,
\emph{On the trivectors
of a }$6$\emph{-dimensional symplectic vector space}, Linear Alg.\
Appl.\ \textbf{435} (2011), 289--306.

\bibitem {BK2} B. De Bruyn, M. Kwiatkowski,
\emph{On the trivectors
of a }$6$\emph{-dimensional symplectic vector space.
II}, Linear Alg.\
Appl.\ \textbf{437} (2012), 1215--1233.

\bibitem {BK3} B. De Bruyn, M. Kwiatkowski,
\emph{On the trivectors
of a }$6$\emph{-dimensional symplectic vector space.
III}, Linear Alg.\
Appl.\ \textbf{438} (2013), 374--398.

\bibitem {BK4} B. De Bruyn, M. Kwiatkowski,
\emph{On the trivectors
of a }$6$\emph{-dimensional symplectic vector space.
IV}, Linear Alg.\
Appl.\ \textbf{438} (2013), 2405--2429.

\bibitem {BK5} B. De Bruyn, M. Kwiatkowski,
\emph{The classification of the trivectors
of a }$6$\emph{-dimensional symplectic space:
Summary, consequences and connections}, Linear Alg.\
Appl.\ \textbf{438} (2013), 3516--3529.

\bibitem {Fogarty} J. Fogarty, \emph{Invariant theory},
W. A. Benjamin, Inc., New York, Amsterdam, 1969.

\bibitem {KPV} V. G. Kac, V. L. Popov,
E. B. Vinberg, \emph{Sur les groupes lin\'eaires
alg\'ebriques dont l'alg\`{e}bre des invariants
est libre}, C. R. Acad. Sci. Paris S\'er. A-B
\textbf{283} (1976), no. 12, Ai, A875--A878.

\bibitem {P} V. L. Popov,
\emph{Classification of spinors of dimension
fourteen}, Trans.\ Mosc.\ Math. Soc.\ \textbf{1}
(1980), 181--232.
\end{thebibliography}
\end{document}